\documentclass[11pt]{article}
\usepackage[utf8]{inputenc}
\usepackage[T1]{fontenc}

\usepackage{amsfonts}
\usepackage{amsmath}
\usepackage{amssymb}
\usepackage{amsthm}

\usepackage{epigraph}

\usepackage{xcolor}
\usepackage{hyperref}

\usepackage{amstext} 
\usepackage{array}   
\newcolumntype{C}{>{$}c<{$}}

\usepackage[shortlabels]{enumitem}

\newtheorem{lem}{Lemma}[section]
\newtheorem{thm}{Theorem}[section]
\newtheorem{cor}{Corollary}[section]
\theoremstyle{definition}
\newtheorem*{df}{Definition}
\newtheorem{conj}{Conjecture}
\newtheorem*{obs}{Observation}

\DeclareMathOperator\Aut{Aut}
\DeclareMathOperator\diam{diam}
\DeclareMathOperator\PSL{PSL}
\DeclareMathOperator\sgn{sgn}
\DeclareMathOperator\ord{ord}

\newcommand{\la}{\langle}
\newcommand{\ra}{\rangle}

\newcommand{\mfr}{\mathfrak}
\newcommand{\mcl}{\mathcal}

\newcommand{\mbb}{\mathbb}
\newcommand{\til}{\widetilde}

\newcommand{\lf}{\left \lfloor}
\newcommand{\rf}{\right \rfloor}

\title{A few remarks on the theory of non-nilpotent graphs}
\author{Radosław Żak\footnote{Undergraduate student;  Jagiellonian University, 
Faculty of Mathematics and Computer Science; 
Łojasiewicza 6, 30-348 Kraków, Poland.\newline \texttt{radoslaw.zak@student.uj.edu.pl}}}
\date{April 22, 2023}

\begin{document}

\maketitle

\begin{abstract}
    We prove a few results about non-nilpotent graphs of symmetric groups $S_n$ -- namely that they satisfy a conjecture of Nongsiang and Saikia (which is likewise proved for alternating groups $A_n$), and that for $n \geqslant 19$ each vertex has degree at least $\frac{n!}{2}$. We also show that the class of non-nilpotent graphs does not have any ``local'' properties, i.e. for every simple graph $X$ there is a group $G$, such that its non-nilpotent graph contains $X$ as an induced subgraph.
\end{abstract}

\section{Introduction}

Let $G$ be a group. We construct a graph $\mcl{N}_G$ with elements of $G$ as vertices, such that $x,\, y \in G$ are connected by an edge if and only if the subgroup $\la x,\, y \ra$ is not nilpotent. By removing isolated vertices of $\mcl{N}_G$ we get the \textit{non-nilpotent graph} of $G$, denoted by $\mfr{R}_G$. Of course in the case when $G$ is weakly nilpotent (i.e. $\la x,\, y \ra$ is nilpotent for any $x,\, y \in G$), $\mfr{R}_G$ is an empty graph. Thus we will be concerned only about non-weakly nilpotent groups; we also will examine only finite groups (for which notions of nilpotency and weak nilpotency coincide).

Non-nilpotent graphs of groups were first defined by A. Abdollahi and M. Zarrin in the paper \cite{abdollahi}, where they proved for instance that in the case of finite non-nilpotent groups $\mfr{R}_G$ are not regular graphs, and that they are connected with diameter bounded by $6$, stating a conjecture that actually $\diam \mfr{R}_G \leqslant 2$. This was disproven by A. Davis, J. Kent and E. McGovern -- in \cite{davis} they have shown that $\diam \mfr{R}_G \geqslant 3$ for $G$ equal to $\mbb{Z}_m \rtimes S_4$, where $m$ is odd and the action of $S_4$ on $\mbb{Z}_m$ is given by $a^\sigma = a^{\sgn \sigma}$. A. Lucchini and D. Nemmi proved in \cite{lucchini} that this bound is indeed optimal, i.e. $\diam \mfr{R}_G \leqslant 3$ for any finite non-nilpotent group $G$.

A natural question is: if $\mfr{R}_G \simeq \mfr{R}_H$ for two finite non-nilpotent groups $G$ and $H$, does it imply that $G \simeq H$? This is not the case -- consider the dihedral group of order $18$, $G = D_9$, and generalized dihedral group of the same order, $H = (\mbb{Z}_3 \times \mbb{Z}_3) \rtimes \mbb{Z}_2$, with action of the non-trivial element of $\mbb{Z}_2$ on $\mbb{Z}_3 \times \mbb{Z}_3$ given by $h \mapsto -h$. In both $G$ and $H$ the only maximal nilpotent subgroups are those cyclic of order $2$ and those of index two (respectively $\mbb{Z}_9$ and $\mbb{Z}_3 \times \mbb{Z}_3$). This means that both $\mfr{R}_G$ and $\mfr{R}_H$ are complete multipartite graphs with nine parts of size one and one part of size eight.

However, a weaker (and hopefully more true) form of this conjecture was stated in \cite{nongsiang} by D. Nongsiang and P. Saikia:

\begin{conj}[Nongsiang \& Saikia] \label{nong}
If $\mfr{R}_G \simeq \mfr{R}_H$ for two finite non-nilpotent groups $G$ and $H$, then $|G|=|H|$.
\end{conj}

They managed to prove this conjecture in several cases, including $G$ being a dihedral group, $\PSL(2,q)$ where $q$ is a prime power with $q^2 \not \equiv 1 \pmod{16}$, groups of order $pq$ and centerless groups of order $pqr$ where $p,\, q,\, r$ are distinct primes. In this article we append this list by two important families of groups: symmetric groups $S_n$ and alternating groups $A_n$ (theorem \ref{SA-theorem}).

It is natural to ask how dense graphs $\mfr{R}_G$ usually are. For a similar notion of the non-commuting graph of $G$, it is easy to see that each vertex is connected with at least half of the other vertices. This is not always true in the non-nilpotent case: $(12)(34)$ has degree $8$ in $\mfr{R}_{S_4}$. However, in the fourth section we prove (Theorem \ref{Deg19}) that graphs $\mfr{R}_{S_n}$ have this property when $n \geqslant 19$. In particular, we can find a Hamiltonian cycle in those graphs by Dirac's theorem.

Much of the attention given to non-nilpotent graphs of groups was focused on cases where $\mfr{R}_G$ has a particularly nice local structure, e.g. if it is a complete multipartite graph. In the last section we prove (theorem \ref{IndSub-Thm}) that we should not hope for these results to be much helpful in the general case: every graph can be found as an induced subgraph of $\mfr{R}_G$ for appropriate $G$.

\section{Notation and preliminaries}

We use the following notations throughout the paper.

\begin{itemize}
    \item When $a$ and $b$ are positive integers, $a \bot b$ means that they are coprime.
    \item The identity element of a group will be usually denoted by $1$.
    \item If $S$ is a subset of a group $G$, then $\la S \ra$ denotes the subgroup generated by $S$.
    \item $G'$ is the commutator subgroup of $G$.
    \item $Z^*(G)$ is the hypercenter (limit of upper central series) of a group $G$.
    \item The semi-direct product of groups $G$ and $H$ is denoted by $G \rtimes_\nu H$ or $G \rtimes H$, where $\nu: H \to \Aut G$ is a homomorphism. (In the second case notation itself does not specify the homomorphism, but we will indicate its choice every time we use it.)
    \item Choose a finite group $G$, an element $g \in G$ and $p$ a prime. Suppose that $g$ has order $n$. Let $n_p$ be equal to $n$ divided by the greatest power of $p$ that divides $n$ -- so that $n_p$ is not divisible by $p$ and $\frac{n}{n_p}$ is a power of $p$. We will denote $g^{n_p}$ by $g_p$. Of course, its order is $\frac{n}{n_p}$, which is a power of $p$. We point out that only finitely many choices of $p$ give $g_p \neq 1$.
    \item If $\sigma$ is a permutation in $S_n$, then its support is the set of elements $x \in \{1,\ldots,n\}$ such that $\sigma(x) \neq x$.

\end{itemize}


We will use the following result extensively:
\begin{lem}[Centralizer structure in symmetric groups] \label{cent-struct}
If $\sigma \in S_n$ has $a_k$ cycles of length $k$, for $k=1,2,\ldots,n$, then the centralizer of $\sigma$ has
\[ \prod_{k=1}^n k^{a_k} a_k! \]
elements, and is a direct product of groups isomorphic to $C_k \wr S_{a_k}$, where $S_{a_k}$ permutes cycles of size $k$ in $\sigma$ (here $C_k$ denotes the cyclic group of order $k$).
\end{lem}

\begin{proof}
Classical. The first part can be found for instance in \cite[p. 78, Ex. 3]{jacobson}; for the second part, we just observe that the direct product is contained in the centralizer and has the same order.
\end{proof}

We assume the following characterizations of finite nilpotent groups to be widely known:

\begin{lem}
Let $G$ be a finite group. The following are equivalent:
\begin{enumerate}
    \item $G$ is nilpotent.
    \item Every Sylow subgroup of $G$ is normal.
    \item $G$ is a direct product of its Sylow subgroups.
    \item If $x,y$ are elements of $G$ with coprime order, then $xy=yx$.
\end{enumerate}
\end{lem}

\begin{proof}
$(1) \iff (2) \iff (3)$ can be found in \cite[p. 24, Thm 1.26]{isaacs}, while $(3) \Rightarrow (4) \Rightarrow (2)$ is straightforward.
\end{proof}

\begin{cor} \label{Pelem2}
Choose $g,h \in G$. Then subgroup $\la g,h \ra$ is nilpotent if and only if the following conditions hold:

\begin{enumerate}
    \item For $p \neq q$ primes, elements $g_p$ and $h_q$ commute.
    \item For any prime $p$, $\la g_p, h_p \ra$ is a $p$-group.
\end{enumerate}
\end{cor}

Let $G$ be a finite group. As in the introduction, we define $\mathcal{N}_G$ to be a graph with vertices being elements of $G$, such that $x$ and $y$ are connected if and only if $\la x,y \ra$ is not nilpotent. The non-nilpotent graph of $G$, denoted by $\mfr{R}_G$, is formed by removing isolated vertices from $\mathcal{N}_G$.

The crucial role of the hypercenter $Z^*(G)$ in investigating $\mfr{R}_G$ was pointed out by Abdollahi and Zarrin in \cite{abdollahi}. They prove in Proposition 2.1 and Lemma 2.2 of their paper the following facts:

\begin{lem}
Isolated vertices of $\mathcal{N}_G$ correspond to the elements of $Z^*(G)$. Moreover, $x$ and $y$ are connected in $\mfr{R}_G$ if and only if their projections on $G/Z^*(G)$ are connected in $\mfr{R}_{G/Z^*(G)}$.
\end{lem}

Quotienting can be rephrased using graph theoretical terms as \textit{blow-ups}. For a graph $X$ and positive integer $c$, denote by $c \odot X$ the graph composed in the following way: for each vertex $a$ of $X$, we put $c$ vertices $a_1,\, a_2,\, \ldots,\, a_c$ in $c \cdot X$, in such a way that $a_i$ and $a_j$ are never connected, and $a_i$ is connected with $b_j$ iff $a$ is connected with $b$ in $X$. This construction is called the blow-up of $X$ of order $c$. The graph $c \odot X$ is the only graph such that there is a $c$-to-$1$ map from $c \odot X$ to $X$ that preserves adjacency relations.

\begin{cor} \label{hyper}
Let $G$ be a finite group. Then $\mfr{R}_G \simeq |Z^*(G)| \odot \mfr{R}_{G/Z^*(G)}$.
\end{cor}

\section{Groups with the same non-nilpotent graph}

Suppose $G$ and $H$ are two finite non-nilpotent groups with $\mfr{R}_G \simeq \mfr{R}_H$. Corollary \ref{hyper} allows us to draw many conclusions from that fact. In particular, by comparing the number of vertices, $|G|-|Z^*(G)|=|H|-|Z^*(H)|$. Therefore instead of proving $|G|=|H|$, we can focus on proving $|Z^*(G)|=|Z^*(H)|$.

\begin{lem} \label{degs}
Let $G$ be a finite non-nilpotent group. Then $|Z^*(G)|$ divides the number of vertices of $\mfr{R}_G$ and the degree of each vertex.
\end{lem}

\begin{proof}
An easy consequence of corollary \ref{hyper}.
\end{proof}

Suppose that the group $G$ is centerless (hence $Z^*(G)$ is trivial) and that $\mfr{R}_G \simeq \mfr{R}_H$ for some finite non-nilpotent group $H$. Then if the greatest common divisor of degrees and the number of vertices in $\mfr{R}_G$ is equal to $1$, we may immediately conclude by lemma \ref{degs} that $|Z^*(H)|=1$, and thus $|G|=|H|$.

It turns out this is precisely the case when $G=S_n$ or $G=A_n$. To prove the Conjecture \ref{nong} in this case, we simply need to look at degrees of certain permutations in the graph $\mfr{R}_G$.

\begin{lem} \label{parts}
Let $\sigma \in G \in \{S_n,A_n\}$. Suppose that $\sigma$ cycle lengths are $p_1$, $p_2$, $\ldots$, $p_k$, where $p_i$ are distinct and either prime or equal to $1$ (hence $\sum p_i = n$). Then for any $\tau \in S_n$ the subgroup $\la \sigma, \tau \ra$ is nilpotent if and only if $\sigma \tau = \tau \sigma$, and there are $\prod p_i$ such permutations $\tau$.
\end{lem}

\begin{proof}
Implication $(\Leftarrow)$ is obvious, so assume $\la \sigma, \tau \ra$ is nilpotent. Choose a prime $q$ and consider $\tau_q$. By Corollary \ref{Pelem2} for all $i$ such that $p_i \neq 1,q$ the element $\sigma_{p_i}$ (which is the corresponding $p_i$-cycle raised to a certain power) commutes with $\tau_q$. By Lemma \ref{cent-struct} $\tau_q$ and $\sigma_{p_i}$ have disjoint supports.

If $q$ is distinct from all $p_i$-s, this immediately forces $\tau_q$ to be trivial. If, however, we have $q=p_i$ for some $i$, then the support of $\tau_q$ is contained in the union of $q$-cycle of $\sigma$ and its possible fixed point. By Corollary \ref{Pelem2} we know that $\la \tau_q, \sigma_q \ra$ is a $q$-subgroup of $S_n$ with support not bigger than $q+1$. But $q^2 \nmid (q+1)!$, so this subgroup is cyclic. Thus $\tau_q$ is generated by $\sigma_q$.

As it is true for all $q$, and they are pairwise coprime, we get $\tau \in \la \sigma \ra$ and the result follows.
\end{proof}

Lemma \ref{parts} provides us with permutations for which we can easily enumerate degrees in $\mfr{R}_G$. A similar family is given by permutations with two cycles of length two and other cycles with lengths being distinct primes, as shown in the following lemma.

\begin{lem} \label{PartTwo}
Let $\sigma \in G \in \{S_n,A_n\}$. Suppose that $\sigma$ cycle lengths are $2$, $2$, $p_1$, $\ldots$, $p_k$, where $p_i$ are distinct odd primes. Then there are $16 \prod p_i$ such permutations $\tau \in S_n$ that $\la \sigma, \tau \ra$ is nilpotent, and $4 \prod p_i$ of them lie in $A_n$.
\end{lem}

\begin{proof}
Without loss of generality, $(12)$ and $(34)$ are the transpositions in $\sigma$.

Similarly to the proof of Lemma \ref{parts}, if $\tau$ is taken in such a way that $\la \sigma, \tau \ra$ is nilpotent, then $\tau_q =1$ for $q \not \in \{2 , p_1, \ldots, p_k\}$, and its support is contained in the support of $\sigma_q$ if $q$ is among these primes, with $\la \tau_q, \sigma_q \ra$ being a $q$-group. For $q=p_i$ we get that $\tau_q$ is a power of the $q$-cycle in $\sigma$. For $q=2$ the element $\tau_2$ acts only on $\{1,2,3,4\}$ in such a way that $\la (12)(34), \tau_2 \ra$ is a $2$-group. We can check that all elements of $S_4$ of order being a power of two can be taken here (due to the three copies of $D_4$ sitting inside $S_4$), and there are $16$ of them (but only four inside $A_4$). In total, we get the desired claim.
\end{proof}

Now we just need to state a final, number-theoretical lemma.

\begin{lem} \label{sum}
Every integer $m \geqslant 28$ can be expressed as a sum of distinct prime numbers greater than $3$.
\end{lem}

\begin{proof}
We can check that all integers in the interval $[28,89]$ can be expressed as a sum of distinct prime numbers from the interval $[5,31]$.

We can proceed by induction. For $n \leqslant 89$, by above, a stronger claim holds, hence we can assume $n \geqslant 90$. By Bertrand's postulate we may find a prime $p$ in the interval $\left[\frac{n-27}{2}, n-27 \right)$. Then we may write $n-p$ as a sum of distinct primes $\geqslant 5$ by induction assumption. It suffices to check that these primes are different from $p$.

As $p \geqslant \frac{n-27}{2} > 31$, we need to worry only about primes generated in the induction steps. But they are decreasing: indeed, in the next step we will choose a prime $q$ in the interval $\left[\frac{n-p-27}{2}, n-p-27 \right)$, and $q < n-p-27 \leqslant \frac{n-27}{2} \leqslant p$.
\end{proof}

\begin{thm} \label{Divs}
Let $n \geqslant 32$ be an integer, and $H$ be a finite group. Then $\mfr{R}_{S_n} \simeq \mfr{R}_H$ implies that $|H|=|S_n|$. Similarly, $\mfr{R}_{A_n} \simeq \mfr{R}_H$ implies that $|H|=|A_n|$.
\end{thm}

\begin{proof}
We just need to prove that $Z^*(H)$ is trivial.

By lemma \ref{sum}, we can write $n-4 = \sum_{i=1}^k p_i$, where $p_i$ are distinct primes greater or equal to $5$. Let $\sigma$ be any permutation with cycle lengths $1$, $3$, $p_1$, $p_2$, $\ldots$, $p_k$. Similarly let $\tau$ be any permutation with cycle lengths $2$, $2$, $p_1$, $p_2$, $\ldots$, $p_k$.

Consider in the first place the case of $S_n$. By Lemma \ref{parts}, $\deg \sigma = n! - 3 \prod p_i$. By Lemma \ref{PartTwo}, $\deg \tau = n! - 16 \prod p_i$. Hence, by Lemma \ref{degs}
\[ |Z^*(H)| \mid \deg \sigma - \deg \tau = 13 \prod p_i. \]
But by the same lemma $|Z^*(H)| \mid n!-1$, so $|Z^*(H)| \bot n!$. However, every prime divisor of $13 \prod p_i$ is also a divisor of $n!$ (as $n \geqslant 13$), thus $|Z^*(H)|=1$.

For $A_n$ using the same results we get $|Z^*(H)| \mid \frac{n!}{2} - 1$, $\deg \sigma = \frac{n!}{2} - 3 \prod p_i$ (note that all permutations mentioned in the Lemma \ref{parts} are even) and $\deg \tau = \frac{n!}{2} - 4 \prod p_i$, hence $|Z^*(H)| \mid \prod p_i$, but again $|Z^*(H)| \bot n!$, thus we again get trivial $Z^*(H)$.
\end{proof}

\begin{thm} \label{SA-theorem}
\begin{itemize}
    \item Let $n \geqslant 3$ be an integer, and $H$ be a finite group. Then $\mfr{R}_{S_n} \simeq \mfr{R}_H$ implies $|H|=|S_n|$.
    \item Let $n \geqslant 4$ be an integer, and $H$ be a finite group. Then $\mfr{R}_{A_n} \simeq \mfr{R}_H$ implies $|H|=|A_n|$.
\end{itemize}
\end{thm}

\begin{proof}
By Theorem \ref{Divs} we just need to check cases $n \leqslant 31$. The only thing we were using in its proof were $n \geqslant 13$ and the ability to express $n-4$ as a sum of distinct primes greater than $3$. This can be done if $n \geqslant 15$ and $n \neq 18,\, 19,\, 25,\, 31$.

Consider the part about $S_n$. Let $z:=|Z^*(H)|$. In general Lemma \ref{parts} means that for any partition of $n=\sum p_i$ as a sum of distinct primes or ones, we have $z \mid n!-\prod p_i$, while Lemma \ref{degs} adds divisibility $z \mid n!-1$. If $n$ has a form $p$ or $p+1$, where $p$ is an odd prime, then those imply $z \mid n!-p$ and $z \mid n!-1$, therefore $z \mid p-1 \mid n!$, thus $z=1$. This leaves us with cases $n=9,\, 10,\, 25$.

For $n=25$ we can use partitions $25=17+7+1=17+5+3$, which imply $z \mid 25! - 17 \cdot 7$, $z \mid 25 ! - 17 \cdot 15$, thus $z \mid 17 \cdot 8$, which is a contradiction with $z \mid 25!-1$. For $n=10$ the partition $10=3+7$ implies $z \mid 10!-21$, and since $z \mid 10!-1$ we get $z \mid 20$, a contradiction. For $n=9$ we take $9=5+3+1$, thus $z \mid 9!-15$ and $z \mid 9!-1$, thus $z \mid 14$, a contradiction.

Note that the above paragraphs work with $A_n$ as well, because all primes we have used in our partitions were odd -- and in this case lemma \ref{parts} works just as well for $A_n$. The only other difference is that we have $\frac{1}{2} n!$ instead of $n!$ in this case, but the reasoning carries through (since $n \geqslant 4$ and $\frac{1}{2}n!$ is even).
\end{proof}

\section{Large degrees}

\begin{thm} \label{Deg19}
Let $n \geqslant 19$ be an integer. Then for every vertex $v$ in $\mfr{R}_{S_n}$ we have $\deg v \geqslant \frac{1}{2} n!$.
\end{thm}

This whole section is aimed to prove Theorem \ref{Deg19}. Some restriction on $n$ being large enough is indeed necessary: we already observed in the previous section that $(12)(34)$ has degree eight in $\mfr{R}_{S_4}$. It is compelling to try shrinking the bound of $19$ using i.e. computer calculations, but it was not attempted by the author.

A direct consequence of Theorem \ref{Deg19} is the existence of Hamiltonian cycles in $\mfr{R}_{S_n}$ for $n \geqslant 19$. Indeed, Dirac's theorem assures us that if a graph has $m$ vertices, each with degree at least $\frac{1}{2}m$, then the graph has a Hamiltonian cycle. Using Bondy-Chvatal theorem \cite{bondychvatal}, which is a generalization of Dirac's theorem, it is possible to prove that $\mfr{R}_{S_n}$ has a Hamiltonian cycle for every $n \geqslant 3$. However, the proof found by the author, while being elementary, is also long and dependent on calculations. Consequently, we decided to omit the proof from this paper; more details can be found in the first version of the article's preprint published on arXiv. The existence of Hamiltonian cycles in non-nilpotent graphs is in line with similar results about non-generating graphs. F. Erdem proved in \cite{erdem} that the non-generating graph of $S_n$ has a Hamiltonian cycle for any $n \geqslant 104$, while earlier Breuer et al. \cite{breuer-etal} obtained the same result for any large enough simple group.

We proceed to the proof of Theorem \ref{Deg19}. Firstly, for any finite group $G$, $x,y \in G$ and $n \in \mathbb{Z}_+$ we have $\la x^n,y \ra \subset \la x,y \ra$, and thus $\deg_{\mathcal{N}_G} x \geqslant \deg_{\mathcal{N}_G} x^n$. Hence

\begin{obs}
It is enough to check Theorem \ref{Deg19} for permutations of prime order.
\end{obs}

Any permutation $\sigma$ of order $p$ is a product of disjoint $p$-cycles, and for any $\tau \in S_n$ we can analyse cycle lengths of $\tau$ to get some results about nilpotence.

\begin{df}
For $\sigma \in S_n$ and $i \in \{1,2,\ldots,n\}$ denote by $c(\sigma,i)$ the length of the cycle in $\sigma$ in which $i$ is contained (in other words, $c(\sigma,i)$ is the smallest $k \in \mathbb{Z}_+$ such that $\sigma^k(i)=i$).
\end{df}

\begin{lem} \label{NilCycle}
Let $\sigma, \tau \in S_n$ be such that $\la \sigma, \tau \ra$ is nilpotent, and $\sigma$ has order $p$ being prime. Take any $i,j \in \{1,2,\ldots,n\}$ that belong to the same cycle in $\sigma$ (with $i \neq j$). Then quotient $\frac{c(\tau,i)}{c(\tau,j)}$ is a rational power of $p$ (so is of form $p^t$, where $t$ is an integer, not necessarily nonnegative).
\end{lem}

\begin{proof}
Consider an element $\rho$ equal to $\tau^{p^k}$, where $p^k$ is the largest power of $p$ dividing $\ord \tau$. Then $\frac{c(\tau,i)}{c(\rho,i)}$ and $\frac{c(\tau,j)}{c(\rho,j)}$ are rational powers of $p$, and the claim will follow if we prove that $c(\rho,i) = c(\rho,j)$. We know that $p \nmid \ord \rho$.

Choose $m$ so that $\sigma^m(i)=j$. Then $\sigma^m$ has order $p$, so it commutes with $\rho$ (as both are elements of nilpotent subgroup $\la \sigma,\, \tau \ra$). Then $c(\rho,i) = c(\sigma^{-m} \rho \sigma^m , i)$, because $\rho = \sigma^{-m} \rho \sigma^m$. But conjugation does not change the length of cycles, while it changes their placement, hence $c(\sigma^{-m} \rho \sigma^m , i) = c(\rho,j)$. Therefore $c(\rho,i) = c(\rho,j)$ and our claim is proven.
\end{proof}

Our strategy now will be to prove that for a fixed $\sigma$ of order $p$, $\tau$ chosen uniformly at random from $S_n$ violates Lemma \ref{NilCycle} with probability at least $\frac{1}{2}$. Therefore we need a result about the distribution of $(c(\tau,i),c(\tau,j))$ for a random $\tau$ and fixed $i,j$. For one coordinate the result is well-known: $c(\tau,i)$ is uniformly distributed in $\{1,2,\ldots,n\}$. Indeed, in exactly $\binom{n-1}{k-1} (k-1)!$ ways we can choose a $k$-cycle containing $i$, so $(n-1)! = \binom{n-1}{k-1} (k-1)!(n-k)!$ permutations have a $k$-cycle containing $i$, and the answer does not depend on $k$.

For two coordinates, the answer still has the same elegant symmetry. We can assume $(i,j)=(1,2)$ without loss of generality.

\begin{lem} \label{rand-perm}
Let $\tau$ be a permutation chosen uniformly at random from $S_n$. Then $1$ and $2$ are in the same cycle with probability $\frac{1}{2}$, and if they are not, then $(c(\sigma,1),c(\sigma,2))$ has a uniform distribution on the set $\{ a,\, b \in \mbb{Z}_+ : a+b \leqslant n \}$.
\end{lem}

\begin{proof}
We start with enumerating the cases where $1$ and $2$ are in the same cycle. If this cycle has length $k \geqslant 2$, then we can choose the other elements of the cycle in $\binom{n-2}{k-2}$ ways, the $k$-cycle on those $k$ elements in $(k-1)!$ ways, and the rest of the permutation in $(n-k)!$ ways. In total, this gives
\[ \frac{(n-2)!}{(k-2)!(n-k)!} \cdot (k-1)! \cdot (n-k)! = (n-2)! (k-1) \]
possibilities. Summing for $k$ from $2$ to $n$ we get
\[ (n-2)! \sum_{k=2}^n (k-1) = (n-2)! \cdot \frac{n(n-1)}{2} = \frac{1}{2} n! \]
options, as claimed.

Now suppose $1$ is in a cycle of length $k$, $2$ is in a cycle of length $l$, and those are distinct cycles (though not necessarily $k \neq l$). As in the previous case, we choose elements of the $k$-cycle in $\binom{n-2}{k-1}$ ways, the cycle itself in $(k-1)!$ ways; the elements of the $l$-cycle in $\binom{n-k-1}{l-1}$ ways, times $(l-1)!$ ways of choosing the cycle, and finally $(n-k-l)!$ for the rest of the permutation. This gives in total
\[\binom{n-2}{k-1} \cdot (k-1)! \cdot \binom{n-k-1}{l-1} \cdot (l-1)! \cdot (n-k-l)! = (n-2)!  \]
options, regardless of the values of $k$ and $l$.
\end{proof}

Lemma \ref{NilCycle} has an additional assumption on $i$ and $j$ lying in the same cycle of $\sigma$. If the same is true for $\tau$ as well, then the lemma becomes useless. Therefore we need some control over that event.

\begin{lem} \label{cykloza}
Let $k \leqslant n$ be positive integers. Each element $\sigma$ of $S_n$ induces a partition on $\{1,\, 2,\ldots,\, k\}$ -- elements in the same cycle of $\sigma$ are put to the same set of the partition. Then, when $\sigma$ is a random permutation in $S_n$, the distribution of partitions will be the same regardless of $n$ (in particular, the same as in the case of $n=k$).
\end{lem}

\begin{proof}
We will proceed by induction on $n$. For $n=k$ there is nothing to prove. Suppose $n>k$ and consider the following function $f: S_n \to S_{n-1}$:
\[ f(\sigma)(j) = \begin{cases} \sigma(j),\quad \sigma(j)\neq n,\\
\sigma(n),\quad \text{else}.
\end{cases} \]
In other words, if we imagined $\sigma$ as a union of disjoint cycles, $f(\sigma)$ is created by deleting $n$ from the permutation, and if this creates a hole in some cycle, we ``sew it up'', putting $\sigma(n)$ as a value of the permutation in the point $\sigma^{-1}(n)$.

It is clear that $f(\sigma)$ yields the same partition on $\{1,\, 2,\ldots,\, k\}$ as $\sigma$ does, thus we only need to prove that $f(\sigma)$ is uniformly distributed whenever $\sigma$ is uniformly distributed.

But $f$ is $n$-to-1 as for any $\til{\sigma} \in S_{n-1}$ and $m \in \{1,\, 2,\ldots,\, n\}$ there is exactly one permutation $\sigma \in S_n$ such that $\til{\sigma} = f(\sigma)$ and $n=\sigma(m)$. For $m=n$ we get it by adding $n$ as a fixed point, and for $m \neq n$ we need to have $\sigma(n)=\til{\sigma}(m)$ and $\sigma$ has to agree with $\til{\sigma}$ outside of $\{m,n\}$.
\end{proof}

Now we have all the machinery needed to prove Theorem \ref{Deg19}.

\begin{proof}[Proof of Theorem \ref{Deg19}.]
We already know it is enough to check that $\deg_{\mathcal{R}_{S_n}} \sigma \geqslant \frac{1}{2} n!$, where $\ord \sigma = p$ is prime. Without loss of generality, $1$ and $2$ are in the same cycle of $\sigma$. We will have several cases.

\begin{enumerate}
    \item $p$ is odd;
    \item $p=2$, $\sigma$ is not a transposition;
    \item $\sigma$ is a transposition.
\end{enumerate}

In each case, we choose $\tau$ uniformly at random from $S_n$, and we claim $\la \sigma, \tau \ra$ is nilpotent with probability at most $\frac{1}{2}$.

\textbf{Case 3.} Suppose $\la \sigma, \tau \ra$ is nilpotent. Let $k$ be an integer such that the order of $\tau^{2^k}$ is odd. Then $\sigma$ has coprime order with $\tau^{2^k}$, hence they commute and by Lemma \ref{cent-struct} we have $c(\tau^{2^k},1)=1$. Therefore $c(\tau,1)$ is a power of two. But $c(\tau,1)$ is uniformly distributed on $\{1,\ldots,n\}$, hence $c(\tau,1)$ is a power of two with probability at most $\frac{\lf \log_2 n \rf + 1}{n}$, which is less than $\frac{1}{2}$ when $n \geqslant 19$ (actually, six is enough).

\textbf{Cases 1. and 2.} For now, assume that $1$ and $2$ lie in different cycles of $\tau$. Then, by Lemma \ref{NilCycle}, $\la \sigma, \tau \ra$ being nilpotent implies that $\frac{c(\tau,1)}{c(\tau,2)}$ is $p^k$ for an integer $k$. Let $(a,b):=(c(\tau,1),c(\tau,2))$. By Lemma \ref{rand-perm} pairs $(a,b)$ are distributed uniformly on the set $S:=\{(a,b) \in \mathbb{Z}_+ : a+b \leqslant n\}$. Set $S$ has $\frac{n(n-1)}{2}$ elements; how many of them have $a/b=p^k$? If $k=0$, we have $a=b \leqslant n/2$, so we have $\lf n/2 \rf$ possibilities. If $k \neq 0$, then we can assume by symmetry $k>0$. Then $a=bp^k$, so $n \geqslant a+b = b(p^k+1)$, which gives us $\lf \frac{n}{p^k+1} \rf$ pairs $(a,b)$.

In total (remembering about cases for negative $k$) we get
\[ h(n,p):= \lf \frac{n}{2} \rf + 2 \sum_{k=1}^\infty \lf \frac{n}{p^k+1} \rf \]
pairs $(a,b) \in S$ such that $a/b$ is a power of $p$. We denote this expression by $h(n,p)$. In this part of the proof, we can deduce that there are at most $h(n,p) \cdot (n-2)!$ permutations $\tau$ such that $\la \sigma, \tau \ra$ is nilpotent and $1$, $2$ are in distinct cycles of $\tau$.

\textbf{Case 1.} In the case when $p$ is odd, without loss of generality we can assume that $1$, $2$ and $3$ are in the same cycle of $\sigma$. Suppose again that $\tau$ is a general permutation such that $\la \sigma, \tau \ra$ is nilpotent.

\begin{itemize}
    \item If $1$ and $2$ are in distinct cycles of $\tau$, we have at most $h(n,p)(n-2)!$ possibilities for $\tau$.
    \item The same happens if $2$ and $3$ are in distinct cycles of $\tau$.
    \item $1$, $2$ and $3$ are in the same cycle for $\frac{1}{3} n!$ different choices of $\tau$, according to Lemma \ref{cykloza}.
\end{itemize}

In total, $\sigma$ is not connected to at most $\frac{n!}{3} + 2h(n,p)(n-2)!$ other vertices of $\mathcal{N}_{S_n}$. We want this to be smaller or equal to $n!/2$, which is equivalent to $h(n,p) \leqslant n(n-1)/12$. But
\[ h(n,p) \leqslant \frac{n}{2} + 2 \sum_{k=1}^\infty \frac{n}{3^k} = \frac{3}{2} n, \]
and this is at most $n(n-1)/12$ as long as $n \geqslant 19$.

\textbf{Case 2.}

We proceed analogously but with different calculations. Without loss of generality, $(12)$ and $(34)$ are two transpositions in $\sigma$. Let $\tau$ be any permutation with $\la \sigma, \tau \ra$ nilpotent.

\begin{itemize}
    \item If $1$ and $2$ are in distinct cycles of $\tau$, we have at most $h(n,2)(n-2)!$ possibilities for $\tau$.
    \item The same number appears if $3$ and $4$ are in distinct cycles of $\tau$.
    \item Assume both $(1,2)$ and $(3,4)$ are in the same cycles of $\tau$. By Lemma \ref{cykloza} there are exactly $\frac{7}{24} n!$ such $\tau$ (because their fraction is the same as in $S_4$, where we have six $4$-cycles plus $(12)(34)$).
\end{itemize}

This in total gives at most $2h(n,2)(n-2)! + \frac{7}{24} n!$ vertices of $\mathcal{N}_{S_n}$ that are not connected to $\sigma$. We want to bound this by $n!/2$ from above, which is equivalent to $h(n,2) \leqslant 5n(n-1)/48$.

However,
\[ h(n,2)/n \leqslant \frac{1}{2} + \frac{2}{3} + \frac{2}{5} + 2\sum_{k=3}^\infty \frac{1}{2^k} = \frac{31}{15}, \]
which is enough for $n \geqslant 21$. For $n=19$ and $20$, we can evaluate $h(n,2)$ directly and check that the claim still holds.

\end{proof}

\section{Induced subgraphs}

It is natural to ask if all graphs $\mfr{R}_G$ share some particular structure. For instance, $\mfr{R}_G$ is a complete multipartite graph for a wide class of cases, including dihedral and dicyclic groups, as well as the groups of order $pq$ with prime $p, q$. The smallest non-example is $\mfr{R}_{S_4}$. 

While being a complete multipartite graph may seem like a global condition, it is actually a local one -- among any three vertices we can find either no or at least two edges. Therefore we can wonder if the class of all graphs $\mfr{R}_G$ can have a characterization by forbidden induced subgraphs, or if we could find any forbidden induced subgraphs at all. This, as the following theorem shows, is not the case.

\begin{thm} \label{IndSub-Thm}
For any simple graph $X$ there is a group $G$ such that $\mfr{R}_G$ has $X$ as an induced subgraph.
\end{thm}

\begin{proof}
Let us label vertices of $X$ by numbers from $\{1,2,\ldots,k\}$. Our construction will go as follows: we will create sequences of groups $\{ 1\} \simeq G_0 < G_1 < \ldots < G_k$ and $N_1,\ N_2,\ \ldots,\ N_k$ so that $G_t = G_{t-1} \ltimes N_t$, and choose $x_i \in N_i\setminus\{1\}$. Our groups $N_t$ will have a form $\mathbb{Z}_{p_t}^{\alpha_t}$, where $p_1$, $p_2$, $\ldots$, $p_k$ are distinct primes to be chosen later. In other words, our main idea is to extend $G$ by an elementary abelian group in each step.

We would want to do it inductively so that the following conditions hold:
\begin{enumerate}
    \item vertices numbered $i$, $j$ are not connected by an edge in $X$ $\iff [x_i,x_j]=1$,
    \item $G_t'$ does not contain any $x_i$ with $1 \leqslant i \leqslant t$.
\end{enumerate}

The first condition is natural and implies that $x_i$-s induce $X$ in $\mfr{R}_G$ as, since their orders are distinct prime numbers, we have $\la x_i,x_j \ra$ nilpotent if and only if $[x_i,x_j]=1$. The second one is a technical assumption important for the construction, but not for the nilpotence itself.

We go by induction on $k$. For $k=1$ choose $N_1=G_1=\mathbb{Z}_2$, and $x_1$ as the non-trivial element.

Now suppose the sequences up to index $t$ are constructed, and we want to perform the inductive step. At first, we consider condition 2. We claim $G_{t+1}' \cap G_t = G_t'$ regardless of how we choose $N_{t+1}$ and extension $G_t \ltimes N_{t+1}$. Direction $\supset$ is easy; to show direction $\subset$ we write epimorphism $\phi: G_{t+1} \to G_t$ (as $G_t \simeq G_{t+1}/N_{t+1}$) and observe that $\phi([x,y]) = [\phi(x),\phi(y)] \in G_t'$. Therefore if $x_r \not \in G_r'$, then $x_r \not \in G_s'$ for any $s \geqslant r$. Thus while performing the step $t \to t+1$ we need to worry only about the case $x_{t+1} \in G_{t+1}'$, which we will sort out at the end.

We claim there is a normal subgroup of $G_t$ (we denote it by $H$) which contains exactly those $x_i$-s, such that vertex $i$ and vertex $t+1$ are not joined by an edge in $X$. Consider the quotient morphism $\mu: G_t \to G_t/G_t'$. We know, by the induction assumption and condition 2, that all $x_i$-s are mapped to non-trivial elements. Group $G_t/G_t'$ is abelian and orders of $x_i$-s (thus also of their images) are distinct primes, therefore we can find a $K \leq G_t/G_t'$ containing images of those $x_i$-s we want -- for instance, any abelian group is a direct product of its Sylow subgroups, and we take subgroups corresponding to the respective primes -- and choose simply $H = \mu^{-1}(K)$. As $K$ is normal in $G_t/G_t'$ (all subgroups of an abelian group are normal), $H$ is normal in $G_t$ by the third isomorphism theorem.

Let $d = [G_t:H]$. We can set $\alpha_{t+1}=d$. By Cayley theorem we have a monomorphism $G_t/H \to S_d$; moreover, we also have a natural embedding $S_d \to \Aut(\mathbb{Z}_{p_{t+1}}^d)$ by permuting the coordinates. Composing these maps together we get
\[ \xi: G_t \to G_t/H \to S_d \to \Aut(\mathbb{Z}_{p_{t+1}}^d) \]
Since only the first arrow is not a monomorphism, $\ker \xi = H$. Now let us define $G_{t+1} = \mathbb{Z}_{p_{t+1}}^d \rtimes_\xi G_t$. Take $x_{t+1}$ to have coordinates that are pairwise different and have a sum different from zero (it can be easily done when $p_{t+1}>d+1$, as either $(1,2,\ldots,d-1,d)$ or $(1,2,\ldots,d-1,d+1)$ will work). Then:

\begin{itemize}
    \item If vertices $i$ and $t+1$ are not connected by an edge in $X$, then $\xi(x_i)$ is trivial as $x_i \in H = \ker \xi$, so $[x_i,x_{t+1}]=1$. 
    \item If they are connected by an edge, then $\xi(x_i)$ is a non-trivial permutation of coordinates on $\mathbb{Z}_{p_{t+1}}^d$, so it acts non-trivially on $x_{t+1}$ (as it has pairwise different coordinates). Therefore $[x_i,x_{t+1}]$ is non-trivial.
\end{itemize}

The only thing left is to ensure that $x_{t+1} \not \in G_{t+1}'$. $G_{t+1}'$ is the intersection of all normal subgroups $K$ of $G_{t+1}$ such that $G_{t+1}/K$ is abelian. Now look at the subset $S$ of $G_{t+1}=N_{t+1} \rtimes G_t$ composed of elements of the form $n * g$ where the sum of coordinates of $n$ is equal to $0$. As, firstly, such $n$-s form a subgroup in $N_{t+1}$ (which is abelian), and secondly, $G_t$ acts without changing the sum of coordinates, hence fixes this set, our set is indeed a normal subgroup. Since $|G_{t+1}/S| = p_{t+1}$, this quotient has to be abelian, thus $G_{t+1}' \subset S \not \ni x_{t+1}$, qed.

This way we have shown that construction of $G_k$ is possible. Now the only thing left is to embed $G_k$ in $G:=S_{|G_k|}$ (by Cayley theorem this is possible), so that all $x_i$-s are outside of $Z^*(G)$ (which is now trivial), so their respective vertices belong not only to $\mathcal{N}_G$, but also to $\mfr{R}_G$ (in case of $k=1$ we need to additionally embed $S_2$ in $S_3$). This ends the proof.
\end{proof}

\section*{Declaration of competing interest}

The authors declare that they have no known competing financial interests or personal relationships that could have appeared to influence the work reported in this paper.

\end{document}